\documentclass[12pt]{amsart}

\usepackage[latin1]{inputenc}
\usepackage[T1]{fontenc}
\usepackage{hyperref}
\usepackage{geometry} 
\usepackage{stmaryrd}
\usepackage{amsmath,amssymb,amsthm}
\usepackage{latexsym}

\newtheorem{theorem}{Theorem}
\newtheorem{proposition}[theorem]{Proposition}
\newtheorem{fact}[theorem]{Fact}
\newtheorem{cor}[theorem]{Corollary}
\newtheorem{lmm}[theorem]{Lemma}
\newtheorem*{claim}{Claim}
\newtheorem{clm}{Claim}

\theoremstyle{definition}
\newtheorem{definition}[theorem]{Definition}
\newtheorem{remark}[theorem]{Remark}
\newtheorem{conj}[theorem]{Conjecture}
\newtheorem{expl}[theorem]{Example}

\def\bsp{\begin{expl}}
\def\ebsp{\end{expl}}
\def\behe{\begin{clm}}
\def\ebehe{\end{clm}}
\def\beh{\begin{claim}}
\def\ebeh{\end{claim}}
\def\defn{\begin{definition}}
\def\edefn{\end{definition}}
\def\satz{\begin{theorem}}
\def\esatz{\end{theorem}}
\def\tats{\begin{fact}}
\def\etats{\end{fact}}
\def\kor{\begin{cor}}
\def\ekor{\end{cor}}
\def\bema{\begin{remark}}
\def\ebema{\end{remark}}
\def\lem{\begin{lmm}}
\def\elem{\end{lmm}}
\def\bem{\begin{remark}}
\def\ebem{\end{remark}}
\def\verm{\begin{conj}}
\def\everm{\end{conj}}
\def\bew{\begin{proof}}
\def\ebew{\end{proof}}
\def\bewbeh{\begin{proof}[Proof of Claim]}
\def\satzli{\begin{proposition}}
\def\esatzli{\end{proposition}}

\def\Ind#1#2{#1\setbox0=\hbox{$#1x$}\kern\wd0\hbox to 
0pt{\hss$#1\mid$\hss}
\lower.9\ht0\hbox to 0pt{\hss$#1\smile$\hss}\kern\wd0}

\def\Notind#1#2{#1\setbox0=\hbox{$#1x$}\kern\wd0\hbox to 0pt{\mathchardef
\nn="3236\hss$#1\nn$\kern1.4\wd0\hss}\hbox to 0pt{\hss$#1\mid$\hss}\lower.9\ht0
\hbox to 0pt{\hss$#1\smile$\hss}\kern\wd0}

\def\indd{\mathop{\ \hbox to 
0pt{$\mid^{\text{d}}$\hss}\,
\lower4pt\hbox to 0pt{\hss$\smile$\hss}\ \ }}
\def\nindd{\mathop{\ \hbox to 
0pt{$\!\not{\mid}^{\text{\,d}}$\hss}\,
\lower4pt\hbox to 0pt{\hss$\smile$\hss}\ 
\ }}

\def\N{\mathbb N}

\def\F{\mathbb F}

\def\Q{\mathbb Q}

\def\semi{\rtimes}
\def\char{\mbox{char}}

\begin{document}

\title{Some remarks on sharply $2$-transitive groups\\ and near-domains}   
\author{Frank O. Wagner}
\date{28/2/2022}

\address{Institut Camille Jordan, Universit\'e Lyon 1}
\email{wagner@math.univ-lyon1.fr}

\begin{abstract}A sharply $2$-transitive permutation group of characteristic $0$ whose point stabiliser has an abelian subgroup of finite index splits. More generally, a near-domain of characteristic $0$ with a multiplicative subgroup of finite index avoiding all multipliers $d_{a,b}$ must be a near-field. In particular this answers question 12.48 b) of the Kourovka Notebook in characteristic 0.\end{abstract}


\maketitle

\section{Introduction}
Recall that a permutation group $G$ acting on a set $X$ is sharply $2$-transitive if for any two pairs $(x,y)$ and $(x',y')$ of distinct elements of $X$ there is a unique $g\in G$ with $gx=x'$ and $gy=y'$. Then $G$ has involutions, and either involutions have fixed points and $G$ is of permutation characteristic $2$, or the action on $X$ is equivalent to the conjugation action on the set $I$ of involutions. In that case, all translations, i.e.\ products of two distinct involutions, are also conjugate and have the same order $p$, which is either an odd  prime or $\infty$; the number $p$ (or $0$ if $p=\infty$) is the permutation characteristic of $G$. We say that $G$ splits if it has a regular normal subgroup $N$; in that case $G=N\semi C_G(i)$ for any involution $i\in I$. Note that Tent, Rips and Segev have constructed non-split sharply $2$-transitive permutation groups of characteristic $0$ and $2$.

V. D. Mazurov asked in the Kourovka Notebook (question 12.48):\\
Let G be a sharply $2$-transitive permutation group.\begin{enumerate}
\item Does G possess a regular normal subgroup if a point stabilizer is locally finite?
\item Does G possess a regular normal subgroup if a point stabilizer has an abelian
subgroup of finite index?\end{enumerate}
We shall answer question (b) affirmatively in permutation charateristic $0$. In fact, we shall show a more general result for near-domains.

\section{Near-domains and near-fields.}
Instead of working with sharply $2$-transitive groups, we shall work in the equivalent setting of near-domains.
\defn
$(K,0,1,+,\cdot)$ is a {\em near-domain} if for all $a,b,c\in K$
\begin{enumerate}
\item $(K,0,+)$ is a {\em loop}, i.e.\ $a+x=b$ and $y+a=b$ have unique solutions, with $a+0=0+a=a$;
\item $(K\setminus\{0\},1,\cdot)$ is a group, and $0\cdot a=a\cdot 0=0$;
\item left distributivity holds: $a\cdot(b+c)=a\cdot b+a\cdot c$;
\item for all $a,b\in K$ there is $d_{a,b}\in K$ such that $a+(b+x)=(a+b)+d_{a,b}\cdot x$ for all $x$.\end{enumerate}
A near-domain is a {\em near-field} if addition is associative.\edefn
Hence a near-field is a skew field iff right distributivity holds.

\tats[Tits, Karzel]
A sharply $2$-transitive permutation group $G$ is isomorphic to the group of affine transformations of some near-domain $K$, i.e.\ of the set of permutations $\{x\mapsto a+bx:a,b\in K,\,b\not=0\}$; the centraliser of any involution is isomorphic to the multiplicative group $K^\times$. It is split iff $K$ is a near-field.\etats

Let $E$ be the set $\{d\in K:1+d=d+1\}$. Since the additive loop of $K$ is power-	associative, it is easy to see that $1$ generates a subfield of $K$ contained in $E$, which is either $\Q$ or $\F_p$. Thus $K$ has a characteristic, which is easily seen to be equal to the permutation characteristic of $G$. Note that in characteristic $>2$ there is a unique maximal sub-near-field, which is equal to $E$.

\tats[\cite{Ker74}] For all $a,b,c\in K$ we have:\begin{enumerate}
\item $d_{a,a}=1$.
\item $d_{a,b}(b+a)=a+b$.
\item $cd_{a,b}c^{-1}=d_{ca,cb}$.
\item $d_{a,b}=d_{a,c}d_{c+a,-c+b}d_{-c,b}$.
\item If $a,b\in E$ then $(a+b)\,2\in E$.
\item $|K^\times:C_{K^\times}(d_{a,b})|=\infty$ if $d_{a,b}\not=1$.
\end{enumerate}
\etats

Let now $A$ be any subgroup of finite index in $K^\times$ which avoids all non-trivial coefficients $d_{a,b}$ for $a,b\in K$.
Kerby \cite[Theorem 8.26]{Ker74} has shown that $K$ must be a near-field in the following cases:\begin{enumerate}
\item $\char K=0$ and $|K^\times:A|=2$,
\item $\char K=2$, $|K^\times:A|=2$ and $|E|>2$,
\item $\char K=p>2$ and $|K^\times:A|<|E|$.
\end{enumerate}
We shall adapt the proof of (3) to characteristic $0$.

\lem\label{lemma} Suppose $d_{a,1/k}=1$. Then $d_{a,n/k}=1$ for all $n\in\N$.\elem
\bew By induction on $n$. This is clear for $n=0$ and $n=1$. So suppose it holds for $n$, and consider
$$\begin{aligned}\frac{n+1}k +a&=\Big(\frac nk+\frac1k\Big)+a=\frac nk+\Big(\frac1k+a\Big)=\frac nk+\Big(a+\frac 1k\Big)\\
&=\Big(\frac nk+a\Big)+\frac 1k=\Big(a+\frac nk\Big)+\frac1k=a+\Big(\frac nk+\frac 1k\Big)=a+\frac{n+1}k.\qedhere\end{aligned}$$
\ebew

\satzli\label{propn} If $A\le K^\times$ is a subgroup of finite index avoiding all nontrivial $d_{a,b}$ and $\char(K)=0$, then $K$ is a near-field.\esatzli
\bew Recall that $\Q\subseteq E$. If $K=E$, then $d_{a,b}=1$ for all $a,b\in K$ and $K$ is a near-field. So assume $E\subsetneq K^\times$, and take $a\in K\setminus E\,2^{-1}$. Let $n=|K^\times:A|$. Then there are distinct $i>j$ in $\{0,1,2,\ldots,n\}/n!$ with $d_{a,i}A=d_{a,j}A$; since $d_{-j,i}=1$ we obtain
$$d_{a,i}=d_{a,j}d_{j+a,-j+i}d_{-j,i}=d_{a,j}d_{j+a,-j+i}.$$
Hence $d_{j+a,-j+i}\in A$, and $d_{j+a,-j+i}=1$ by assumption.

Now $d_{(i-j)^{-1}(j+a),1}=d_{j+a,-j+i}=1$, so $(i-j)^{-1}(j+a)\in E$. Since $-(i-j)^{-1}j\in\Q\subseteq E$, we have 
$$[-(i-j)^{-1}j+(i-j)^{-1}(j+a)]\,2=(i-j)^{-1}a\,2\in E,$$
and $d_{a2,i-j}=1$. But $0<(i-j)\,n!\le n$ is integer, and there is an integer $k>0$ with $i-j=\frac1k$. By Lemma \ref{lemma} we obtain $d_{a2,1}=1$ and $a\,2\in E$, a contradiction.
\ebew

\kor Let G be a sharply doubly transitive permutation group of characteristic $0$ whose point stabilizer is virtually abelian. Then $G$ is split.\ekor
\bew If $K$ is the associated near-domain, $K^\times$ has an abelian subgroup $A$ of finite index. Now any non-trivial $d_{a,b}$ has a centralizer of infinite index in $K^\times$, so $d_{a,b}\notin A$. We finish by Proposition \ref{propn}.\ebew


\begin{thebibliography}{99}

\bibitem{Ker74} Kerby, William. {\em On infinite sharply multiply transitive groups}, Hamburger Mathematische Einzelschriften, Neue Folge, Heft 6. Vandenhoeck \& Ruprecht, G\"ottingen 1974.
\end{thebibliography}
\end{document}